\documentclass[runningheads]{llncs}
\usepackage[T1]{fontenc}
\usepackage{amsmath}
\usepackage{amssymb}
\usepackage[hidelinks]{hyperref}
\usepackage{cleveref}
\usepackage{listings}
\usepackage{algorithm}
\usepackage{algpseudocode}

\crefname{lstlisting}{listing}{listings}
\Crefname{lstlisting}{Listing}{Listings}
\AtBeginDocument{}
\let\origthelstnumber\thelstnumber
\makeatletter
\newcommand*\Suppressnumber{%
  \lst@AddToHook{OnNewLine}{%
    \let\thelstnumber\relax%
  }%
}

\newcommand*\Reactivatenumber{%
  \lst@AddToHook{OnNewLine}{%
   \let\thelstnumber\origthelstnumber%
  }%
}
\makeatother
\usepackage{graphicx}
\usepackage{siunitx}
\usepackage{subcaption}
\usepackage{xcolor}
\usepackage{xspace}
\usepackage[mode=tex,subpreambles=true]{standalone}
\usepackage{color}

\mathchardef\mhyphen="2D
\newcommand{\flexi}{FLEXI\xspace}
\newcommand{\pyhope}{PyHOPE\xspace}
\newcommand{\hopr}{HOPR\xspace}

\begin{document}
%
\title{In-Memory Load Balancing for Discontinuous Galerkin Methods on Polytopal Meshes}
\titlerunning{Load Balancing for DG Methods on Polytopal Meshes}
%
\author{Patrick Kopper\orcidID{0000-0002-7613-0739} \and
Anna Schwarz\orcidID{0000-0002-3181-8230} \and
Jens Keim\orcidID{0000-0002-2338-1497} \and
Andrea Beck\orcidID{0000-0003-3634-7447}}
\authorrunning{P. Kopper et al.}
%
\institute{Institute of Aerodynamics and Gas Dynamics, University of Stuttgart, Stuttgart, Germany \\
\email{\{kopper,schwarz,keim,beck\}@iag.uni-stuttgart.de}}
\maketitle              
\begin{abstract}
High-order accurate discontinuous Galerkin (DG) methods have emerged as powerful tools for solving partial differential equations such as the compressible Navier--Stokes equations due to their excellent dispersion-dissipation properties and scalability on modern hardware.
The open-source DG framework \flexi has recently been extended to support DG schemes on general polytopal elements including tetrahedra, prisms, and pyramids.
This advancement enables simulations on complex geometries where purely hexahedral meshes are difficult or impossible to generate.
However, the use of meshes with heterogeneous element types introduces a workload imbalance, a consequence of the temporal evolution of modal rather than nodal degrees of freedom and the accompanying transformations.
In this work, we present a lightweight, system-agnostic in-memory load balancing strategy designed for high-order DG solvers.
The method employs high-precision runtime measurements and efficient data redistribution to dynamically reassign mesh elements along a space-filling curve.
We demonstrate the effectiveness of the approach through simulations of the Taylor--Green vortex and large-scale parallel runs on the EuroHPC pre-exascale system MareNostrum~5.
Results show that the proposed strategy recovers a significant fraction of the lost efficiency on heterogeneous meshes while retaining excellent strong and weak scaling.

\end{abstract}
\section{Introduction}
Compressible, turbulent flows pose a significant challenge towards numerical simulations due to their wide range of flow scales and the accompanying resolution requirements to resolve intricate flow details.
High-order methods alleviate some of these challenges due to their low points-per-wavenumber requirements and favorable dispersion-dissipation properties, maximizing utilization of the amount of available random access memory and the often constrained interconnect bandwidth.
One particularly efficient high-order scheme is the discontinuous Galerkin spectral element method (DGSEM) which uses an
element-local weak form of the underlying flow equations to formulate performant algorithms on modern hardware.
Traditionally, DGSEM utilized quadrilateral and hexahedral elements in order to apply tensor-product operators in a dimension-by-dimension manner, resembling the highly efficient approach of spectral methods~\cite{Orszag1979}.
Although DGSEM yields a computationally attractive scheme, the automatic generation of purely hexahedral meshes for complex geometries remains an area of ongoing research~\cite{Blacker2001,Shepherd2008}.
As such, numerous authors have explored the effective extension of DGSEM to non-hexahedral
possibly curved elements~\cite{Chan2016,Montoya2024a} by the use of a collapsed coordinate transformation based on a \textit{Duffy transformation}~\cite{Duffy1982}.
In addition, the authors were able to circumvent the restrictive time step limitation frequently encountered in collapsed coordinate
discretizations~\cite{Dubiner1991} by means of the temporal evolution of modal rather than nodal degrees of freedom, called modal
time stepping.
While these approaches reduce the computational complexity of non-hexahedral elements from $\mathcal{O}(p^{2d})$ to
$\mathcal{O}(p^{d+1})$, non-hexahedral elements still introduce a workload imbalance through the requirement of additional matrix
multiplications for the modal time stepping. 
In this paper, we present an easy and generally applicable workload balancing approach for a high-order DG solver based on space-filling Hilbert curves (SFC).
SFC-based approaches limit the time spend in identifying a new workload distribution by reducing the complex topological locality
problem posed by a three-dimensional unstructured grid to a 1D space-partitioning problem and have been successfully applied to
several high-order solvers~\cite{Harlacher2012,Borrell2018,Reinarz2020,Kopper2023,Nayak2025}.
This approach is highly efficient through the use of in-memory workload balancing, yet remains system-agnostic by using high-precision runtime timers.
The approach is implemented in the open-source solver framework \flexi\footnote{\url{https://github.com/flexi-framework/flexi}}~\cite{Krais2019}, developed by the Numerics Research Group\footnote{\url{https://numericsresearchgroup.org}} located at the Institute of Aerodynamics and Gas Dynamics at the University of Stuttgart.

\section{Numerical Methods}
\flexi solves the compressible unsteady Navier--Stokes--Fourier equations via high-order DGSEM, which is particularly well suited for high-performance computing as it combines the geometric flexibility of unstructured meshes with computationally efficient tensor-product operations.
The method exhibits high arithmetic intensity and excellent cache reuse, rendering it attractive for modern many-core architectures.
In DGSEM, the domain is tessellated into non-overlapping, possibly non-conforming and potentially curved elements using unstructured, potentially moving grids.
The governing equations are transformed from the physical space into the polytopal reference coordinate system
$\boldsymbol{\mathcal{T}} \subset \mathbb{R}^3$ of the underlying element type via the mapping $\mathbf{x} = \boldsymbol{\chi}(\boldsymbol{\xi})$.
The variational form is then obtained via $L_2$ projection of the governing equations onto the space of element-local Legendre--Gauss polynomial test functions up to degree $\mathcal{N}$.
The weak form of the governing equations is obtained after application of Gauss's theorem and reads
\begin{align}
  \int_{         \boldsymbol{\mathcal{T}} }{J \frac{\partial\mathbf{q}_{h}}{\partial t}}\phi({\boldsymbol{\xi}}) \,d{\boldsymbol{\xi}} +
  \int_{\partial{\boldsymbol{\mathcal{T}}}}(\boldsymbol{\mathcal{F}}\cdot\mathbf{n})^{*}\phi({\boldsymbol{\xi}})           \,d{\boldsymbol{S}}
 -\int_{         \boldsymbol{\mathcal{T}} }{\mathcal{F}}(\mathbf{q}_{h},\nabla\mathbf{q}_{h})\cdot\nabla_{\xi}\phi({\boldsymbol{\xi}}) \,d{\boldsymbol{\xi}} = 0,
\end{align}
where $\boldsymbol{q}_h = \boldsymbol{q}_h(\boldsymbol{\xi},t)$ is the element-local solution, $\boldsymbol{\mathcal{F}}$ the contravariant flux vector, $J$ the Jacobian of the mapping $\boldsymbol{\chi}$, and $\boldsymbol{n}$ the outward-pointing normal vector.
For non-hexahedral elements, the polytopal reference space is transformed again to the unit hexahedral via the Duffy transformation~\cite{Duffy1982}.
Further details on the derivation of the DGSEM on polytopal elements can be found in~\cite{Keim2025}.
The DGSEM is integrated in time using low-storage explicit Runge--Kutta (LSERK) schemes following the method of lines approach.
To alleviate the common drawback of DG methods, as is the case with other high-order methods, their lack of robustness in the
presence of strongly non-linear flux functions (aliasing), it is necessary to employ additional stabilization techniques.
For this, an entropy-stable formulation of the DGSEM based on a summation-by-parts like property and an adequate two-point
flux~\cite{Chandrashekar2013} is utilized~\cite{Schwarz2025,Keim2025}.
\flexi is written in modern Fortran and parallelized using pure MPI following the MPI-everywhere paradigm using non-blocking communication and the MPI-3 shared memory (SHM) scheme.

\paragraph{Modal Time Stepping}
Explicit time integration of the nodal DGSEM on polytopal element types leads to prohibitory small time steps.
This issue is circumvented by advancing the modal polynomial coefficients in time, enabling stable and accurate integration at larger time steps.
For this modal time stepping, the nodal degrees of freedom are mapped into modal space using a generalized Vandermonde matrix, integrated in time using the same LSERK approach as the nodal coefficients, and the solution transformed back to nodal space using the inverse Vandermonde.
Again, see~\cite{Keim2025,Montoya2024,Montoya2024a} and \Cref{alg:modal} for details.
However, the modal time stepping comes at the cost of increased computational overhead due to the additional Vandermonde
transformations, introducing a workload imbalance across varying element types.
\begin{algorithm}
  \caption{Schematic procedure for modal time stepping\label{alg:modal}}%
  \begin{algorithmic}
    \State{Precompute Vandermonde matrix $V = \phi(\boldsymbol{\xi})$, physical mass matrix $M = V^T W J V$ and weight matrix $W$}
    \While{$t < t_\mathrm{end}$}
      \State{Solve a projection problem to obtain modal representation $\tilde{q}_h = M^{-1} V^T W J q_h$ (cannot be
      precomputed since the Jacobian can change over time)}
      \State{Evolve modal degrees of freedom in time $\tilde{q}_h \gets$ LSERK($\tilde{q}_h$)}
      \State{Compute nodal degrees of freedom $q_h = V \tilde{q}_h$}
      \State{Update $t=t+\Delta t$}
    \EndWhile
  \end{algorithmic}
\end{algorithm}

\paragraph{Software Infrastructure}
\flexi is maintained with modern software engineering practices to ensure correctness and reproducibility.
A continuous integration (CI) pipeline automatically builds the code and runs a suite of regression tests\footnote{\url{https://github.com/reggie-framework/reggie2.0}} whenever changes are pushed to the repository.
These tests cover a broad range of element types, polynomial orders, and parallel configurations, guaranteeing that new features or optimizations do not introduce regressions or performance degradations.

\section{Parallelization Strategy}
\begin{figure}[!tb]
  \begin{subfigure}{0.6375\linewidth}
  \includegraphics[width=\linewidth]{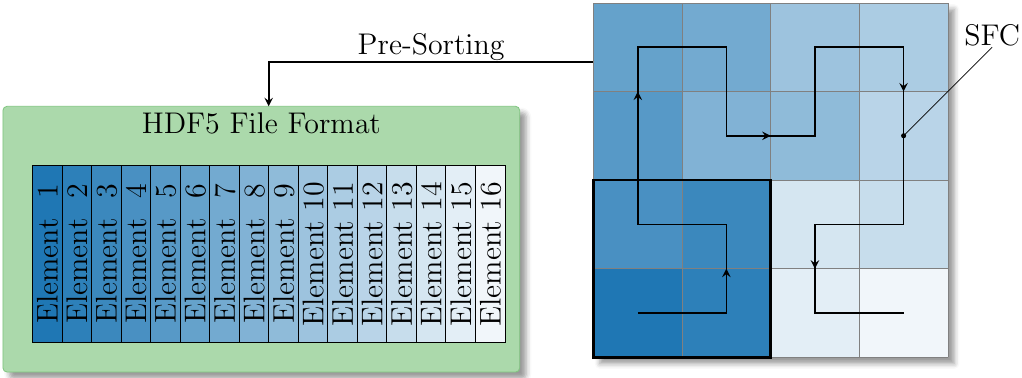}
    \caption{Mesh pre-sorting}%
    \label{fig:mpi:sfc}
  \end{subfigure}\hfill
  \begin{subfigure}{.3125\linewidth}
  \includegraphics[width=\linewidth]{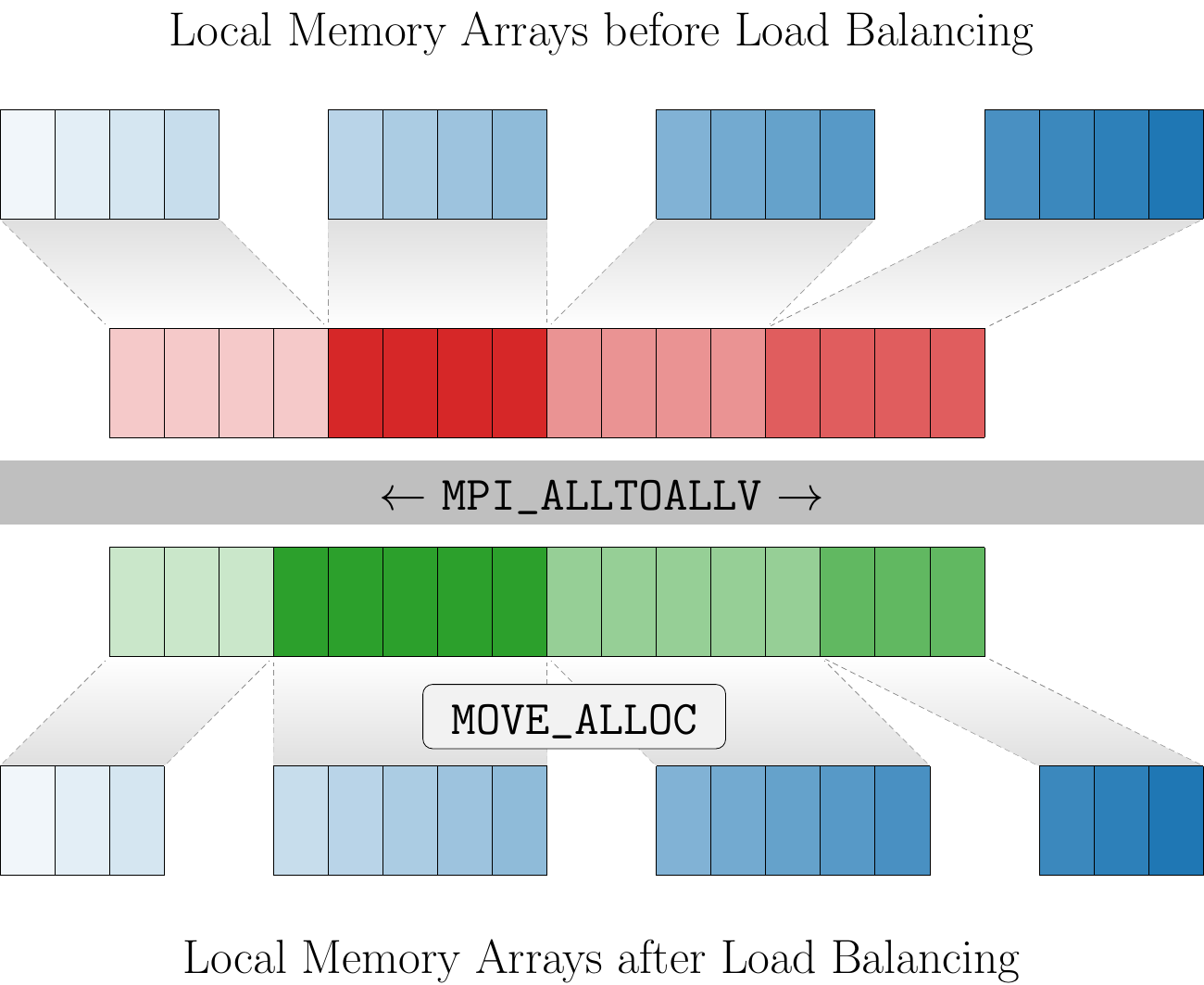}
    \caption{Data exchange}%
    \label{fig:mpi:exchange}
  \end{subfigure}
  \caption{Mesh pre-sorting along the space-filling curve and in-memory data exchange during load balancing.}%
\end{figure}

The element-local basis of the DGSEM leads to a straightforward parallelization approach as the coupling between elements is restricted to the exchange of the numerical flux.
\flexi relies on its mesh preprocessor \pyhope\footnote{\url{https://pypi.org/project/PyHOPE}}~\cite{Kopper2025} to sort the unstructured grid elements along a space-filling Hilbert curve (SFC) during mesh generation.
Mesh connectivity information together with face and node information is stored non-uniquely along the same SFC in the \hopr HDF5 file format~\cite{Hindenlang2014}, see \Cref{fig:mpi:sfc}.
Compared to graph-based partitioning approaches, SFC strategies rely only on geometric locality.
While graph-based methods can minimize edge cuts by incorporating connectivity information, they often incur significant preprocessing costs and require global communication.
During code execution in \flexi, the SFC-sorted grid is partitioned into continuous, non-overlapping segments of the SFC which are assigned to individual MPI ranks and loaded using non-overlapping data access.
Geometric information is furthermore stored uniquely on a node-level using SHM arrays, see~\cite{Kopper2022}.

The modal time stepping of the polytopal element types leads to a significant load imbalance as the transformation from nodal to modal space and inverse can increase the computational workload of an individual element by up to \SI[round-mode=places, round-precision=0]{46.20644610319181}{\percent}.
To counter this workload imbalance, \flexi is equipped with a purely in-memory restart-based load balancing approach.
Conditional high-precision timers are invoked at fixed intervals to measure the actual computational effort of each operation loop, cf.\ \Cref{listing:timer}.
\begin{lstlisting}[language=fortran,
                   morekeywords={MPI_WTIME,SideToElem},
                   frame=lines,
                   float=!tb,
                   escapeinside=||,
                   caption=Fortran example code for element- and side-based timers,
                   label=listing:timer,
                  ]
tStart=MPI_WTIME()
CALL VolumeIntegral()                                    ! Element-based
t(DG_ELEMS)=t(DG_ELEMS)+MPI_WTIME()-tStart; tStart=MPI_WTIME()
CALL SurfaceIntegral()                                   ! Side-based
t(DG_SIDES)=t(DG_SIDES)+MPI_WTIME()-tStart; tStart=MPI_WTIME()
CALL ModalTimestep()                                     ! Element-based
t(DG_MODAL)=t(DG_MODAL)+MPI_WTIME()-tStart|\Suppressnumber|
|\ldots||\Reactivatenumber|
DO iElem=1,nElems
  tElem(iElem)=tElem(iElem)+t(DG_ELEMS)/nElems
  IF(ModalElem(iElem)) tElem(iElem)=tElem(iElem)+t(DG_MODAL)/nModalElems
END DO
DO iSide=1,nSides
  iElem=SideToElem(iSide)
  tElem(iElem)=tElem(iElem)+t(DG_SIDES)/nSides
END DO
\end{lstlisting}
The recorded runtime is assigned to individual elements based on their loop participation and load balancing is triggered once the load imbalance exceeds an acceptable threshold.
As the speed-up gained from performing the load balance must offset the time spent with re-initialization, this threshold is inherently system-dependent.
The actual load exchange is performed via shifting of elements along the SFC using collective \verb!MPI_ALLTOALLV! calls to reduce the amount of data sent, followed by a \verb!MOVE_ALLOC! call to move the allocation to the new distribution without invoking an additional memory copy, see \Cref{fig:mpi:exchange,listing:mpi}.
\begin{lstlisting}[language=fortran,
                   morekeywords={ASSOCIATE,ENDASSOCIATE},
                   frame=lines,
                   float=!tb,
                   caption=Fortran example code for load balancing of the solution array \texttt{U},
                   label=listing:mpi,
                  ]
ALLOCATE(UTmp(nVar,0:N,0:N,0:N,nElems))
ASSOCIATE(&
  countSend => (nVar*(N+1)*(N+1)*(N+1)*MPInElemSend     ),&
  dispSend  => (nVar*(N+1)*(N+1)*(N+1)*MPIoffsetElemSend),&
  countRecv => (nVar*(N+1)*(N+1)*(N+1)*MPInElemRecv     ),&
  dispRecv  => (nVar*(N+1)*(N+1)*(N+1)*MPIoffsetElemRecv)
  CALL MPI_ALLTOALLV(U   ,countSend,dispSend,MPI_DOUBLE_PRECISION,&
                     UTmp,countRecv,dispRecv,MPI_DOUBLE_PRECISION,&
                     MPI_COMM,iError)
END ASSOCIATE
CALL MOVE_ALLOC(UTmp,U)
\end{lstlisting}
It is important to note that while an SFC-based decomposition approach exhibits superior speed compared to graph-based methods, these methods lack the ability to incorporate connectivity information when determining a new distribution.

\section{Parallel Performance}
To evaluate the effectiveness of the runtime load balancing strategy in \flexi, we consider both single-node and large-scale distributed memory experiments.
A detailed discussion of the simulation accuracy of \flexi using varying element types is given in~\cite{Keim2025}.
Two representative test cases are investigated: the three-dimensional Taylor--Green vortex, which serves as a controlled benchmark for analyzing workload imbalance on purely hexahedral and heterogeneous element type meshes, and a large-scale advection problem on the EuroHPC pre-exascale system MareNostrum~5 to assess strong and weak scaling behavior.

\begin{figure}[!tb]
  \centering%
  \begin{subfigure}{.325\linewidth}
    \includegraphics[width=\linewidth]{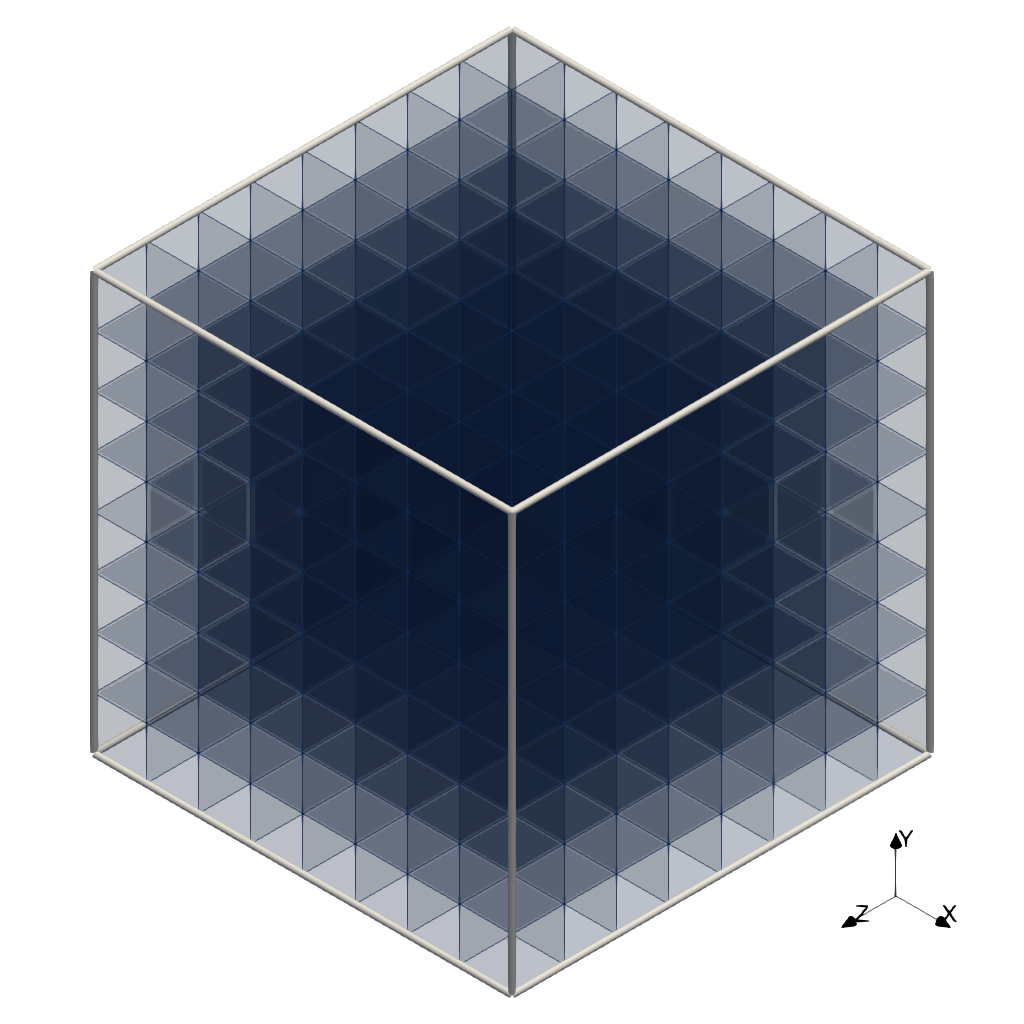}
    \caption{Computational grid}
  \end{subfigure}\hfill
  \begin{subfigure}{.65\linewidth}
    \includegraphics[width=\linewidth]{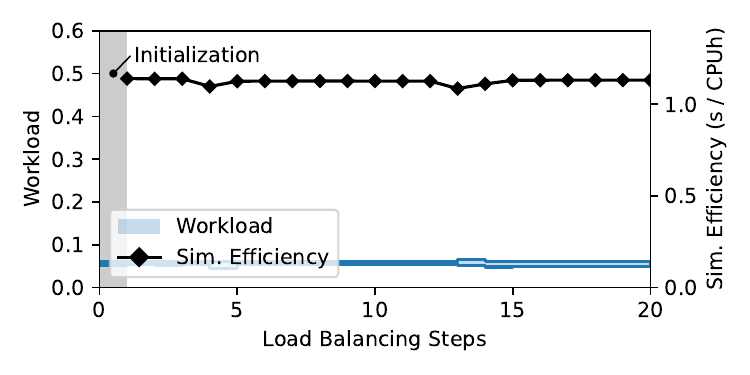}
    \caption{Workload imbalance and performance index}
  \end{subfigure}
  \caption{Computational grid and parallel performance for the Taylor--Green vortex on \num{512} element mesh with purely hexahedral elements.}%
  \label{fig:tgv:hexa}
\end{figure}

\begin{figure}[!tb]
  \centering%
  \begin{subfigure}{.325\linewidth}
    \includegraphics[width=\linewidth]{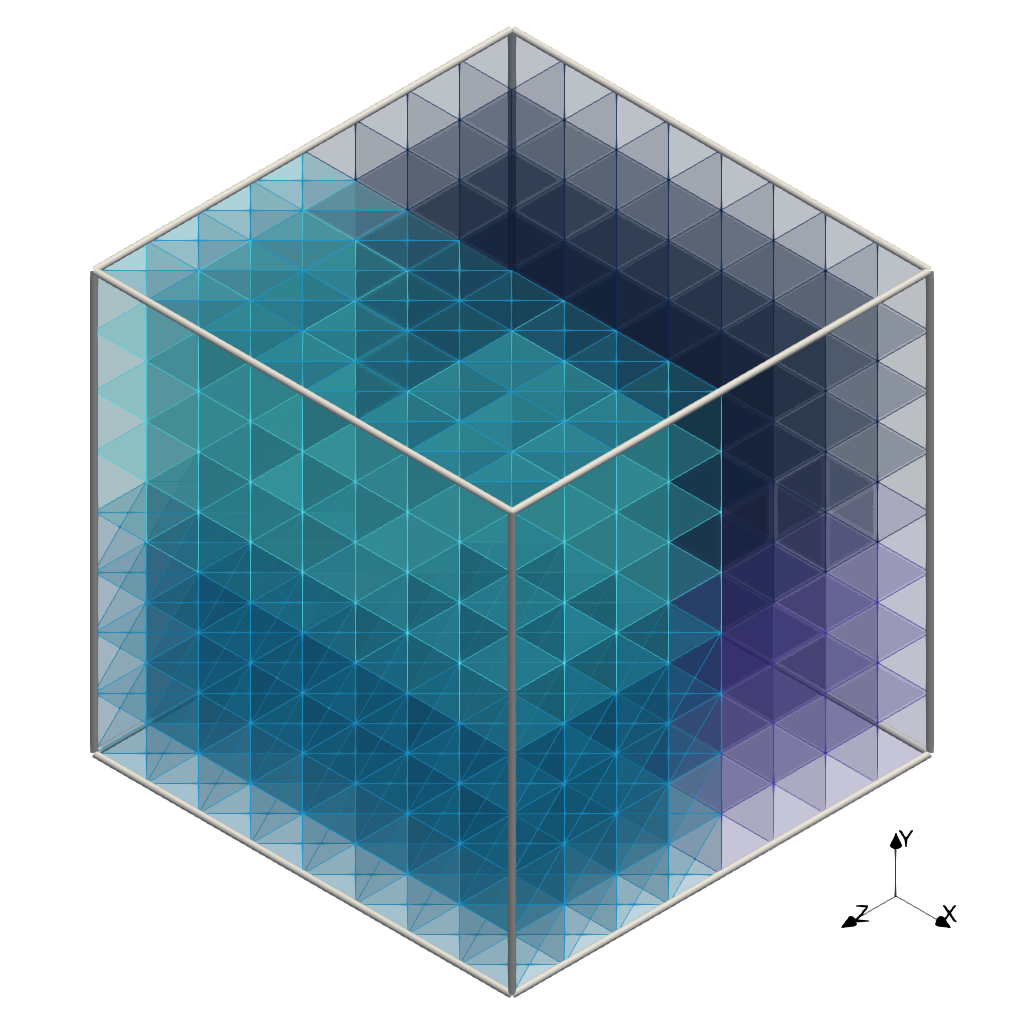}
    \caption{Computational grid}
  \end{subfigure}\hfill
  \begin{subfigure}{.65\linewidth}
    \includegraphics[width=\linewidth]{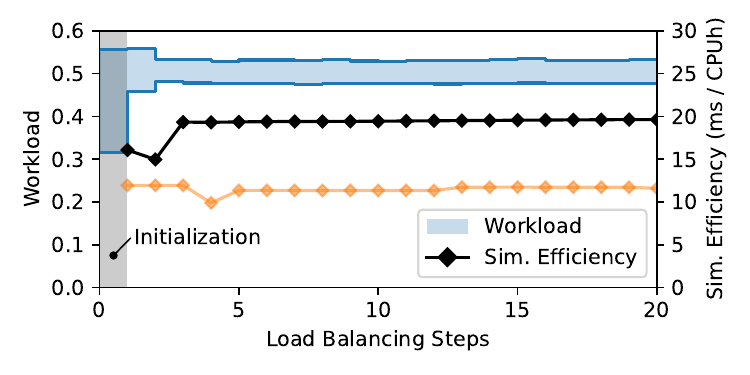}
    \caption{Workload imbalance and performance index}
  \end{subfigure}
  \caption{Computational grid and parallel performance for the Taylor--Green vortex on \num{1984} element mesh with mixed element types. The element types are distinguished by their color. The baseline without load balancing is shaded in orange.}%
  \label{fig:tgv:mixed}
\end{figure}

\paragraph{Taylor--Green Vortex}
The parallel performance of \flexi is evaluated first via simulation of the three-dimensional, viscous Taylor--Green vortex (TGV) flow at Reynolds number $Re = \num{1600}$ on a single compute node using a dual-socket Intel Xeon Gold 6140 CPU configuration with \num{18} cores per socket and \SI{96}{\giga\byte} of RAM.
The TGV was originally proposed by Taylor and Green~\cite{Taylor1937} to study the turbulent energy cascade and features a temporally decaying turbulent flow field undergoing transition and relaxation inside a periodic cube.
Parallel performance is evaluated using both the workload difference in terms of non-dimensional recorded runtime per MPI rank and the simulation efficiency.
The latter represents the simulated time advanced by the explicit Runge--Kutta time stepping scheme for each hour of CPU time spent.
The TGV on the purely hexahedral grid with \num[scientific-notation=false, round-mode=places, round-precision=0]{512} elements, shown in \Cref{fig:tgv:hexa}, results in an even load distribution with a simulation efficiency of \SI[round-mode=places, round-precision=1]{1.0821923752927123}{\second\per{CPU}\hour} for $\mathcal{N} = 5$.
\Cref{fig:tgv:mixed} displays the workload imbalance and simulation efficiency for the TGV on an \num{1984} element mesh with mixed element types.
The grid features the same dimensions as the purely hexahedral case but \num{384} of the \num{512} the hexahedral elements were split to introduce \num{896} tetrahedrons, \num{704} pyramids, and \num{256} prisms or wedges.
The resultant increase in the total number of mesh elements inherently contributes to a decrease in overall simulation efficiency due to the larger problem size.
Furthermore, while the workload imbalance is evident, workload balancing is not trivially successful on distributed systems due to the inability to take communication starvation into account when redistributing along the SFC.
Nonetheless, the load balancing strategy is partly able to recover over time when compared to the baseline simulation without load balancing (shaded in orange).
Note that the simulation efficiency of \SI[round-mode=places, round-precision=1]{18.0986198728929}{\micro\second\per{CPU}\hour} for
$\mathcal{N} = 5$ stays below the value for the purely hexahedral grid, with the difference partly stemming from the required modal time stepping and partly attributed to the varying number of mesh elements within each grid.

\begin{figure}[!tb]
  \includegraphics[width=\linewidth]{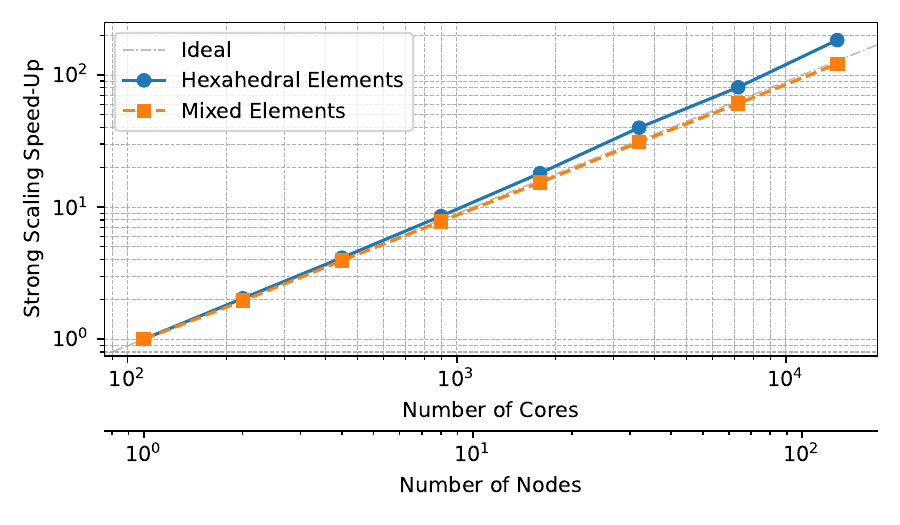}
  \caption{Strong scaling speed-up for the mesh with purely hexahedral elements and the mesh with mixed element types.}%
  \label{fig:scaling:strong}
\end{figure}

\begin{figure}[!tb]
  \includegraphics[width=\linewidth]{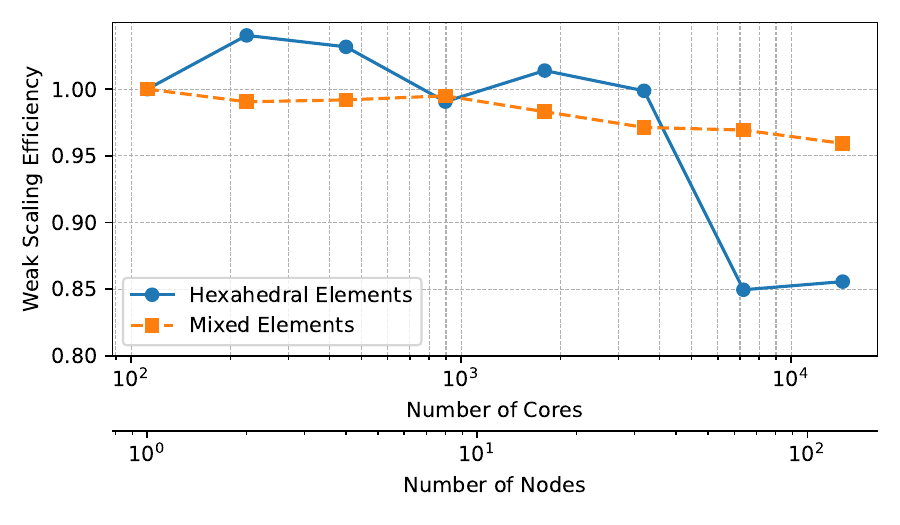}
  \caption{Weak scaling efficiency for the mesh with purely hexahedral elements and the mesh with mixed element types.}%
  \label{fig:scaling:weak}
\end{figure}

\paragraph{Large-Scale Performance}
The large-scale performance of \flexi is evaluated on the pre-exascale EuroHPC supercomputer MareNostrum 5 hosted at BSC-CNS.
MareNostrum 5 features a dual-socket Intel Xeon Platinum 8480+ configuration with \num{56} cores per socket and \SI{256}{\giga\byte} of RAM per node, connected with InfiniBand NDR200 in a fat-tree topology.
Simulations were performed by calculating the advection of an undisturbed flow state using a grid with elemental dimension $L \times L \times L$.
Mixed element type meshes were generated from the purely hexahedral mesh by splitting each quarter of the cross-section into an
individual element type.
As a metric for performance, the performance index (PID) is employed, defined as
\begin{equation}
	\mathrm{PID} = \frac{\mathrm{wall}\mhyphen\mathrm{clock}\mhyphen\mathrm{time} \cdot \mathrm{\#ranks}}{\mathrm{\#DOF} \cdot \mathrm{\#time\;steps} \cdot \mathrm{\#RK}\mhyphen\mathrm{stages}},
\end{equation}
and is a measure of the average time taken by a computational rank to update a single degree of freedom for one Runge--Kutta
stage.
In all simulations, the PID was normalized to the single-node value and averaged over \num{5} runs to eliminate machine and interconnect fluctuations.
\Cref{fig:scaling:strong} shows the strong scaling performance of \flexi for meshes with \num{262144} hexahedral elements and \num{495616} elements of mixed type, respectively.
Both cases demonstrate excellent strong scaling with the purely hexahedral grid even exhibiting superlinear scaling due to reduced cache pressure.
Weak scaling tests were performed by subsequently doubling the grid dimension and the number of elements in streamwise direction.
Weak scaling efficiency on MareNostrum 5 is depicted in \Cref{fig:scaling:weak} and shows a more differentiated behavior.
While the purely hexahedral grid exhibits similar superlinear scaling as for the strong scaling case for \num{32} nodes and below, the weak scaling efficiency drops sharply to about \SI{85}{\percent} for \num{64} nodes and above as communication starvation becomes the dominating factor.
Running \flexi on meshes with mixed element types retains a high weak scaling efficiency above \SI{95}{\percent} throughout all simulation cases considered as the higher single-node load serves to hide communication latency at high node counts.

\section{Conclusion}
In this work, a lightweight, in-memory load balancing strategy for high-order DG solvers operating on meshes with heterogeneous element types is presented.
The approach leverages high-precision runtime measurements and redistributes mesh elements along a space-filling curve to mitigate the workload imbalance introduced by modal time stepping and Vandermonde transformations.
Performance studies on both single-node and large-scale distributed memory systems demonstrate that the method effectively reduces imbalance in heterogeneous meshes while maintaining the strong and weak scaling characteristics of the underlying DG solver.
Although the load balancing does not fully recover the efficiency of purely hexahedral meshes, it provides a practical and system-agnostic means to improve performance in complex geometries without significant memory or communication overhead.
Overall, the proposed strategy extends the applicability of \flexi to large-scale simulations on mixed-element meshes, providing a flexible framework for high-order computations on modern high-performance computing architectures.

\begin{credits}
\subsubsection{\ackname}
The research presented in this paper was funded in parts by Deutsche Forschungsgemeinschaft (DFG, German Research
Foundation) under Germany's Excellence Strategy - EXC 2075 (390740016) and by the European Union and by the state of Baden-Württemberg under the project Aerospace 2050 MWK32-7531-49/13/7 "FLUTTER".
This work has received funding from the European High Performance Computing Joint Undertaking (JU) and Sweden, Germany, Spain, Greece, and Denmark under grant
agreement No 101093393.
We acknowledge the support by the Stuttgart Center for Simulation Science (SimTech).
We acknowledge the EuroHPC Joint Undertaking for awarding us access to MareNostrum5 at BSC, Spain.

\subsubsection{\discintname}
The authors have no competing interests to declare that are
relevant to the content of this article.
\end{credits}
%
%
\bibliographystyle{splncs04}
\bibliography{references}

@article{Blacker2001,
  author        = {Blacker, T.},
  journal       = {Engineering With Computers},
  title         = {Automated Conformal Hexahedral Meshing Constraints, Challenges and Opportunities},
  year          = 2001,
  issn          = {0177-0667},
  month         = oct,
  number        = 3,
  pages         = {201--210},
  volume        = 17,
  doi           = {10.1007/pl00013384},
  publisher     = {Springer Science and Business Media LLC},
}

@article{Borrell2018,
  author        = {Borrell, R. and Cajas, J. C. and Mira, D. and Taha, A. and Koric, S. and V{\'{a}}zquez, M. and Houzeaux, G.},
  journal       = {Computers \& Fluids},
  title         = {Parallel mesh partitioning based on space filling curves},
  year          = 2018,
  month         = sep,
  pages         = {264--272},
  volume        = 173,
  doi           = {10.1016/j.compfluid.2018.01.040},
  publisher     = {Elsevier {BV}},
}

@article{Chan2016,
  author        = {Chan, Jesse and Wang, Zheng and Modave, Axel and Remacle, Jean-Francois and Warburton, T},
  doi           = {10.1016/j.jcp.2016.04.003},
  journal       = {Journal of Computational Physics},
  pages         = {142--168},
  publisher     = {Elsevier Inc.},
  title         = {{GPU-accelerated discontinuous Galerkin methods on hybrid meshes}},
  volume        = 318,
  year          = 2016,
}

@article{Chandrashekar2013,
  author        = {Chandrashekar, Praveen},
  journal       = {Communications in Computational Physics},
  title         = {Kinetic Energy Preserving and Entropy Stable Finite Volume Schemes for Compressible {Euler} and {Navier}-{Stokes} Equations},
  year          = 2013,
  issn          = {1991-7120},
  month         = nov,
  number        = 5,
  pages         = {1252--1286},
  volume        = 14,
  doi           = {10.4208/cicp.170712.010313a},
  publisher     = {Global Science Press},
}

@article{Dubiner1991,
  author        = {Dubiner, Moshe},
  doi           = {10.1007/BF01060030},
  journal       = {Journal of Scientific Computing},
  number        = 4,
  pages         = {345--390},
  title         = {{Spectral methods on triangles and other domains}},
  volume        = 6,
  year          = 1991,
}

@article{Duffy1982,
  author        = {Duffy, Michael G.},
  doi           = {10.1137/0719090},
  journal       = {SIAM Journal on Numerical Analysis},
  number        = 6,
  pages         = {1260--1262},
  title         = {{Quadrature Over a Pyramid or Cube of Integrands with a Singularity at a Vertex}},
  volume        = 19,
  year          = 1982,
}

@inproceedings{Harlacher2012,
  author        = {Harlacher, Daniel F. and Klimach, Harald and Roller, Sabine and Siebert, Christian and Wolf, Felix},
  booktitle     = {2012 {IEEE} 26th International Parallel and Distributed Processing Symposium Workshops {\&} {PhD} Forum},
  title         = {Dynamic Load Balancing for Unstructured Meshes on Space-Filling Curves},
  year          = 2012,
  month         = may,
  publisher     = {IEEE},
  doi           = {10.1109/ipdpsw.2012.207},
}

@phdthesis{Hindenlang2014,
  author        = {Hindenlang, Florian},
  title         = {Mesh curving techniques for high order parallel simulations on unstructured meshes},
  year          = 2014,
  doi           = {10.18419/opus-3957},
  publisher     = {University of Stuttgart},
}

@article{Keim2025,
  title         = {Entropy stable high-order discontinuous {G}alerkin spectral-element methods on curvilinear, hybrid meshes},
  author        = {Keim, Jens and Schwarz, Anna and Kopper, Patrick and Blind, Marcel and Rohde, Christian and Beck, Andrea},
  journal       = {arXiv preprint arXiv:2507.04334},
  year          = 2025,
}

@article{Kopper2022,
  author        = {Kopper, Patrick and Copplestone, Stephen M. and Pfeiffer, Marcel and Koch, Christian and Fasoulas, Stefanos and Beck, Andrea},
  journal       = {Advances in Engineering Software},
  title         = {Hybrid parallelization of {E}uler{\textendash}{L}agrange simulations based on {MPI}-3 shared memory},
  year          = 2022,
  month         = dec,
  pages         = 103291,
  volume        = 174,
  doi           = {10.1016/j.advengsoft.2022.103291},
  publisher     = {Elsevier {BV}},
}

@article{Kopper2023,
  author        = {Kopper, Patrick and Schwarz, Anna and Copplestone, Stephen M. and Ortwein, Philip and Staudacher, Stephan and Beck, Andrea},
  journal       = {Computer Physics Communications},
  title         = {A framework for high-fidelity particle tracking on massively parallel systems},
  year          = 2023,
  month         = apr,
  pages         = 108762,
  volume        = 289,
  doi           = {10.1016/j.cpc.2023.108762},
  publisher     = {Elsevier {BV}},
}

@misc{Kopper2025,
  author        = {Kopper, Patrick and Blind, Marcel P. and Schwarz, Anna and Kurz, Marius and Rodach, Felix and Copplestone, Stephen M. and Beck, Andrea D.},
  title         = {{PyHOPE}: A {P}ython Toolkit for Three-Dimensional Unstructured High-Order Meshes},
  note          = {Manuscript submitted for publication},
  year          = 2025,
}

@article{Krais2019,
  title         = {{FLEXI}: A high order discontinuous {G}alerkin framework for hyperbolic-parabolic conservation laws},
  author        = {Krais, Nico and Beck, Andrea and Bolemann, Thomas and Frank, Hannes and Flad, David and Gassner, Gregor and Hindenlang, Florian and Hoffmann, Malte and Kuhn, Thomas and Sonntag, Matthias and Munz, Claus-Dieter},
  year          = 2021,
  journal       = {Computers {\&} Mathematics with Applications},
  volume        = 81,
  pages         = {186--219},
  doi           = {10.1016/j.camwa.2020.05.004},
  issn          = {08981221},
  archiveprefix = {arXiv},
  arxivid       = {1910.02858},
  eprint        = {1910.02858},
}

@article{Montoya2024,
  author        = {Montoya, Tristan and Zingg, David W.},
  journal       = {Journal of Computational Physics},
  title         = {Efficient entropy-stable discontinuous spectral-element methods using tensor-product summation-by-parts operators on triangles and tetrahedra},
  year          = 2024,
  issn          = {0021-9991},
  month         = nov,
  pages         = 113360,
  volume        = 516,
  doi           = {10.1016/j.jcp.2024.113360},
  publisher     = {Elsevier BV},
}

@article{Montoya2024a,
  author        = {Montoya, Tristan and Zingg, David W.},
  journal       = {SIAM Journal on Scientific Computing},
  title         = {Efficient Tensor-Product Spectral-Element Operators with the Summation-by-Parts Property on Curved Triangles and Tetrahedra},
  year          = 2024,
  issn          = {1095-7197},
  month         = jul,
  number        = 4,
  pages         = {A2270--A2297},
  volume        = 46,
  doi           = {10.1137/23m1573963},
  publisher     = {Society for Industrial \& Applied Mathematics (SIAM)},
}

@article{Nayak2025,
  author        = {Nayak, Amit and Mavriplis, Catherine},
  journal       = {Computers \&amp; Fluids},
  title         = {Immersed boundaries in the discontinuous {G}alerkin spectral element method through hp-adaptivity},
  year          = 2025,
  issn          = {0045-7930},
  month         = nov,
  pages         = 106840,
  volume        = 302,
  doi           = {10.1016/j.compfluid.2025.106840},
  publisher     = {Elsevier BV},
}

@inbook{Orszag1979,
  author        = {Orszag, Steven A.},
  pages         = {273--305},
  publisher     = {Elsevier},
  title         = {Spectral methods for problems in complex geometrics},
  year          = 1979,
  isbn          = 9780125460507,
  booktitle     = {Numerical Methods for Partial Differential Equations},
  doi           = {10.1016/b978-0-12-546050-7.50014-9},
}

@article{Reinarz2020,
  author        = {Reinarz, Anne and Charrier, Dominic E. and Bader, Michael and Bovard, Luke and Dumbser, Michael and Duru, Kenneth and Fambri, Francesco and Gabriel, Alice-Agnes and Gallard, Jean-Matthieu and K\"{o}ppel, Sven and Krenz, Lukas and Rannabauer, Leonhard and Rezzolla, Luciano and Samfass, Philipp and Tavelli, Maurizio and Weinzierl, Tobias},
  journal       = {Computer Physics Communications},
  title         = {ExaHyPE: An engine for parallel dynamically adaptive simulations of wave problems},
  year          = 2020,
  issn          = {0010-4655},
  month         = sep,
  pages         = 107251,
  volume        = 254,
  doi           = {10.1016/j.cpc.2020.107251},
  publisher     = {Elsevier BV},
}

@article{Shepherd2008,
  author        = {Shepherd, Jason F. and Johnson, Chris R.},
  journal       = {Engineering with Computers},
  title         = {Hexahedral mesh generation constraints},
  year          = 2008,
  issn          = {1435-5663},
  month         = mar,
  number        = 3,
  pages         = {195--213},
  volume        = 24,
  doi           = {10.1007/s00366-008-0091-4},
  publisher     = {Springer Science and Business Media LLC},
}

@article{Taylor1937,
  author        = {Taylor, Geoffrey Ingram and Green, Albert Edward},
  journal       = {Proceedings of the Royal Society of London. Series A - Mathematical and Physical Sciences},
  title         = {Mechanism of the production of small eddies from large ones},
  year          = 1937,
  month         = feb,
  number        = 895,
  pages         = {499--521},
  volume        = 158,
  doi           = {10.1098/rspa.1937.0036},
  publisher     = {The Royal Society},
}

@article{Schwarz2025,
author = {Schwarz, Anna and Kempf, Daniel and Keim, Jens and Kopper, Patrick and Rohde, Christian and Beck, Andrea},
title = {{Comparison of Entropy Stable Collocation High-Order DG Methods for Compressible Turbulent Flows}},
url = {http://arxiv.org/abs/2504.00173},
journal = {Computers {\&} Fluids},
year = {2025}
}
\end{document}